\documentclass[11pt]{article}

\usepackage{graphicx,subfigure}
\usepackage{epsfig}
\usepackage{amssymb}
\usepackage{mathrsfs}
\usepackage{epstopdf}
\usepackage{color}
\usepackage{amsmath}
\usepackage{ifpdf}

\newcommand{\E}{{\cal E}}

\newtheorem{theorem}{Theorem}[section]
\newtheorem{definition}{Definition}[section]
\newtheorem{lemma}[theorem]{Lemma}

\def\whitebox{{\hbox{\hskip 1pt
 \vrule height 6pt depth 1.5pt
 \lower 1.5pt\vbox to 7.5pt{\hrule width
    3.2pt\vfill\hrule width 3.2pt}%
 \vrule height 6pt depth 1.5pt
 \hskip 1pt } }}
\def\qed{\ifhmode\allowbreak\else\nobreak\fi\hfill\quad\nobreak
     \whitebox\medbreak}

\newcommand{\ignore}[1]{}

\begin {document}

\baselineskip 16pt
\title{Bicyclic graphs with extremal degree resistance distance}

 \author{\small   Jia\textrm{-}Bao \ Liu$^{a,b}$, \ Si\textrm{-}Qi Zhang$^{b}$, \ \ Xiang\textrm{-}Feng \ Pan$^{b}$\thanks{Corresponding author. Tel:+86-551-63861313.
  \E-mail:liujiabaoad@163.com(J.-B.Liu), zhangsiqidb@163.com(S.-Q. Zhang), xfpan@ahu.edu.cn(X.-F. Pan), shaohuiwang@yahoo.com, swang4@go.olemiss.edu(S. Wang), sakander1566@gmail.com(S. Hayat).},  \ \ Shaohui \ Wang$^c$, \ \ Sakander Hayat$^{d}$\\
 \small  $^{a}$ School of Mathematics and Physics,
Anhui Jianzhu University, Hefei, 230601, P.R. China\\
 \small  $^{b}$ School of Mathematical Sciences, Anhui  University, Hefei 230601, P.R. China\\
 \small  $^{c}$ Department of Mathematics, The University of Mississippi, University, MS 38677, USA\\
 \small $^{d}$ School of Mathematical Sciences, University of Science and Technology of China (USTC), \\ \small Hefei, 230026, P.R. China}

\date{}
\maketitle
\begin{abstract}
Let $r(u,v)$ be the resistance distance between two vertices $u,
v$ of a simple graph $G$, which is the effective resistance
between the vertices in the corresponding electrical network
constructed from $G$ by replacing each edge of $G$ with a unit
resistor. The degree resistance distance of a simple graph $G$ is
defined as ${D_R}(G) = \sum\limits_{\{ u,v\}  \subseteq V(G)}
{[d(u) + d(v)]r(u,v)},$ where $d(u)$ is the degree of the vertex
$u$. In this paper, the bicyclic graphs with extremal degree
resistance distance are strong-minded. We first determine the
$n$-vertex bicyclic graphs having precisely two cycles with minimum
and maximum degree resistance distance. We then completely
characterize the bicyclic graphs with extremal degree resistance
distance.

\end{abstract}



\section{ Introduction}

Topological indices based on the various distances between the
vertices of a graph are widely used in theoretical chemistry to
establish relations between the structure and the properties of
molecules. They provide correlations with physical, chemical and
thermodynamic parameters of chemical
compounds~\cite{Gao2012,Li2015,Tomescu2015,XuG2014,Yang2014,YangJ2008,YangK2013,YangZhang2008}.

The graphs considered in this paper are finite, loopless, and
containing no multiple edges. Given a graph $G$, let $V(G)$ and
$E(G)$ be, respectively, its vertex and edge sets. The ordinary
distance $d(u,v)= d_G(u,v)$ between the vertices $u$ and $v$ of
the graph $G$ is the length of the shortest path between $u$ and
$v$~\cite{Bondy1976}.  The Wiener index~\cite{Gao2011} of a
connected graph $G$, denoted by $W(G)$, is defined as the sum of
all distances between unordered pairs of vertices $u$ and $v$,
i.e.,
$$W(G)=\sum_{\{u,v\}\subseteq V(G)}d(u,v).$$
The resistance distance between the vertices $u$ and $v$ of the
graph $G$, denoted by $r(u,v)$, is defined to be the effective
resistance between the nodes $u$ and $v$ in $G$. Analogous to
Wiener index, the Kirchhoff index is defined
as~\cite{DSTG2015,Klein1993}
\[Kf(G)=\sum\limits_{\{ u,v\} \subseteq V(G)} {r(u,v)}
\]which has been widely studied~\cite{Bu2014,AFeng2012,FengY2014,Feng2014,FengL2014,Feng2012,
LY2015,Liu2014,LP2015,LiuPL2015,LiuYan,Palacios2004,Palacios2015,Palacios2010,Palacios2001,Palacios2011,PalaciosJ2010,Hua2010,Xu2015,YangK2013,YangK2015,ZZ}.
 For other notations in graph theory, we
follow~\cite{Bondy1976} and recent
papers~\cite{LiuC2015,LiuCA2015,Liu2015,LPA2015,LiuP2015,wang2015}.

 A modified version of the
Wiener index is the degree distance defined as
$$D(G)=\sum_{\{u,v\}\subseteq V(G)}[d(u)+d(v)]d(u,v),$$ where
$d(u)= d_G(u)$ is the degree of the vertex $u$ of the graph $G.$
It was introduced independently by Gutman~\cite{Gutman1996},
Dobrynin and Kochetova~\cite{AFeng2012} as a weighted version of
the Wiener index.


Analogous to the definition of the degree distance, the degree
resistance distance is defined as
\[{D_R}(G) = \sum\limits_{\{ u,v\}  \subseteq
V(G)} {[d(u) + d(v)]r(u,v)},\] which was introduced by I. Gutman,
L. Feng and G. Yu in~\cite{Feng2012}.
 Palacios~\cite{Palacios2013} named
the same graph invariant ``additive degree-Kirchhoff index" and
gave tight upper and lower bounds for the degree resistance
distance of a connected undirected graph by using Markov chain
theory.

Tomescu~\cite{Tomescu2008} determined the unicyclic and bicyclic
graphs with minimum degree distance $D(G)$. In~\cite{Tomescu2009}, the
author investigated the properties of connected graphs having
minimum degree distance $D(G)$.
 Bianchi et al.~\cite{Bianchi2013} gave some upper and lower bounds for degree
resistance distance $D_R(G)$ whose expressions do not depend on
the resistance distances. Recently, Chen et al.~\cite{chensubo}
characterized the maximal degree resistance distance of unicyclic
graphs. In a natural way, we will consider the bicyclic graphs. A
bicyclic graph is a connected graph whose edge number is one more
than its vertex number. Obviously a bicyclic graph contains either
two or three cycles. Throughout this article, we restrict our
consideration on bicyclic graphs with exactly two cycles.

For convenience, let $\mathscr {B}_n^{p,q}$ be the set of the
bicyclic graphs with exactly two cycles $C_p$ and $C_q$ as
follows: $C_p=v_1v_2\cdots v_pv_1$ and $C_q=u_1u_2\cdots u_qu_1$
are two cycles such that there is a path $P=v_1w_1\cdots
w_{m-1}u_1$ joining them. The trees $T_{v_i}$, $T_{u_j}$ and
$T_{w_k}$ are rooted at $v_i$, $u_j$ and $w_k$, respectively.
 We say a tree $T$ trivial if $|V(T)|=1$, i.e., $T$ is an singleton vertex, see Fig. 1.

\begin{minipage}[t]{0.5\linewidth}
  \centering
    \includegraphics[scale=0.75]{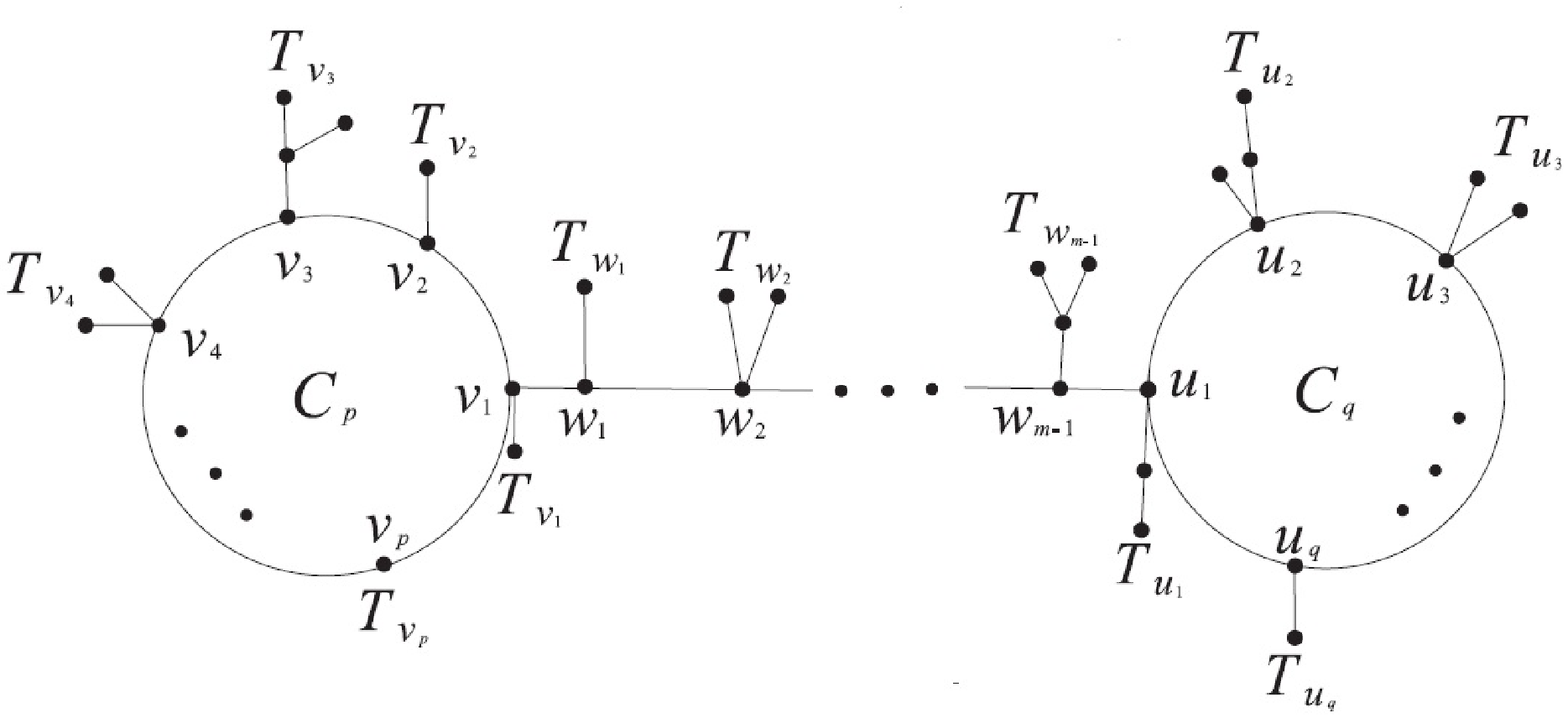}
   \centerline{ ~~~~~~~~~~~~~~~~~~~~~~~~~~~~~~~~~~~~~~~~~~~~~~~~~~~~~~~~~Fig. 1. Illustration for a graph in graph $\mathscr {B}_n^{p,q}$.}
\end{minipage}%
\vspace{2em}

 Let $S^{p,q}_n$ be the graph obtained from cycles $C_p$ and $C_q$ by attaching $n+1-p-q$ pendent edges to the unique common vertex of them.
 Let $P^{p,q}_n$ be the graph consisting of two disjoint cycles $C_p$ and $C_q$ and a path of length $n-p-q+1$ joining them.
 If $n-p-q+1=0$, $P^{p,q}_n$ coincides with $S^{p,q}_n$, see Fig. 2.

\begin{minipage}[t]{0.5\linewidth}
  \centering
    \includegraphics[scale=0.75]{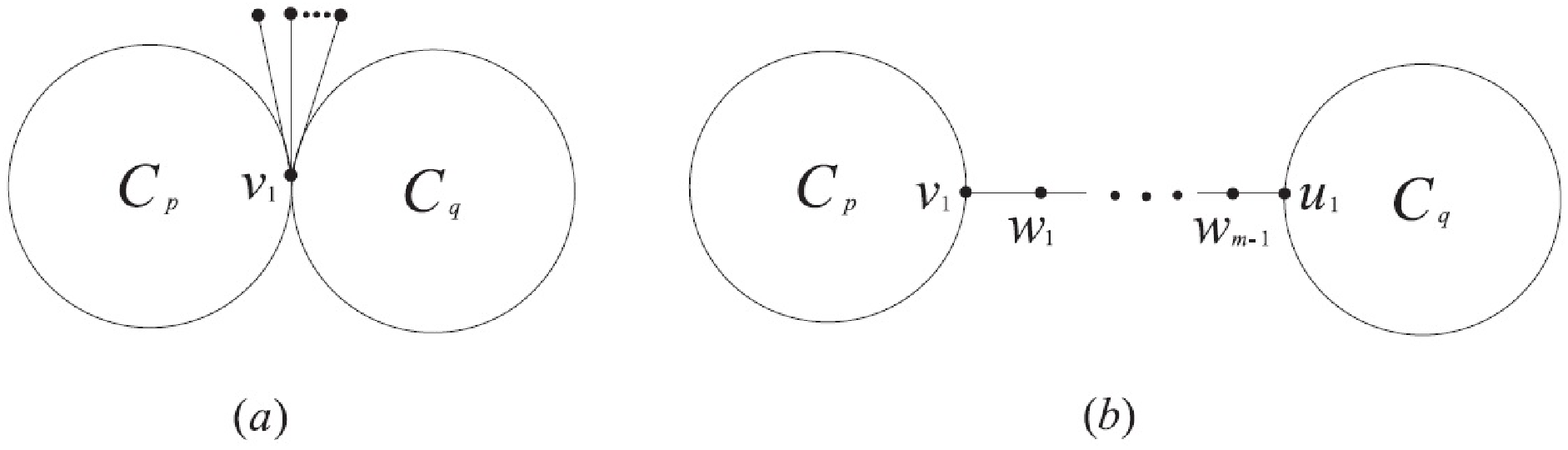}
   \centerline{ ~~~~~~~~~~~~~~~~~~~~~~~~~~~~~~~~~~~~~~~~~~~~~~~~~~~~~~~~~Fig. 2. (a) $S^{p,q}_n$  and (b)
$P^{p,q}_n$.}
\end{minipage}%
\vspace{2em}

To the best of our knowledge, the extremal degree resistance
distances for bicyclic graphs has not been considered so far. In
this paper, we firstly characterize $n$-vertex bicyclic graphs
 with exactly two cycles having minimum and maximum degree resistance
distance, and then characterize the bicyclic graphs with extremal
degree resistance distance.

\section{Preliminaries}

In this section, we provide some lemmas and graph transformations,
that play an important role in the subsequent sections.
 Let ${r}(u,v)$ denote the resistance distance between $u$ and
$v$ in the graph $G$. Recall that ${r}(u,v) = {r}(v,u)$ and
${r}(u,v) \ge 0$ with equality if and only if $u = v$.

In what follows, for the sake of simplicity, instead of $u \in
V(G)$ we write $u \in G $. For a vertex $v$ in $G$, we define
\begin{center}
  $K{f_v}(G) = \sum\limits_{u \in G} {{r}(u,v)}$ ~and~ ${D_v}(G) = \sum\limits_{u \in G} {{d}(u){r}(u,v)}$.
\end{center}
The definition of ${D_v}(G)$ implies that \[{D_R}(G) =
\sum\limits_{v \in G} {{d}(v)} \sum\limits_{u \in G} {{r}(u,v).}
\]

For a vertex $v\in G$, $G-v$ denotes the graph obtained from $G$
by deleting $v$ and its incident edges. Let $\overline{G}$ be the
complement of $G$. Let $G-e$ ($G+e$) denote the graphs obtained
from $G$ by deleting (or adding respectively) the edge $e$.

  Let $H$ be a subgraph of graph $G$. For a vertex $u\in
  V(H)$, let$$r(u|H)=\sum\limits_{v \in V(H)} r({v,u|H}),~ S^{'}(u|H)=\sum\limits_{v \in V(H)}
  {d(v)r(v,u|H)}.$$

  If $C_n$ is a cycle with the vertex set $V(C_n) = \{v_1,\dots,v_n\}$ and $n\geq3$, by Ohm's law, we
  have that, for
  $v_i,v_j\in V(C_n)$ with $i<j$,$$r_{C_n}(v_i,v_j)=\frac{(j-i)(n+i-j)}{n}.$$

\begin{lemma}(\cite{Klein1993}).
 Let $x$ be a
cut-vertex of a connected graph, and let $a$ and $b$ be the
vertices occurring in different components which arise upon
deletion of $x$. Then $$ r(a,b)=r(a,x)+r(x,b).$$
\end{lemma}
\begin{lemma} (\cite{Feng2012}).
 Let $C_k$ be the cycle of size k,
 and $v \in C_k$.
 Then

 $$Kf({C_k}) = \frac{{{k^3} - k}}{{12}},~
 {D_R}({C_k}) = \frac{{{k^3} - k}}{3},~
 K{f_v}({C_k}) = \frac{{{k^2} - 1}}{6},~ {D_v}({C_k}) = \frac{{{k^2} - 1}}{3}.$$
\end{lemma}

\begin{lemma} (\cite{Feng2012}).
  Let $G_1$ and $G_2$ be connected graphs with disjoint vertex sets, with $n_1$ and $n_2$ vertices,
  and with $m_1$ and $m_2$ edges, respectively. Let $u_1\in V(G_1), ~u_2\in V(G_2)$.
  Constructing the graph $G$ by identifying the vertices $u_1$ and $u_2$, and denote the so obtained vertex by $u$.
  Then \[{D_R}(G) = {D_R}({G_1}) + {D_R}({G_2}) + 2{m_2}r(u_1|{G_1}) + 2{m_1}r(u_2|{G_2})
   + ({n_2} - 1)S^{'}(u_1|{G_1}) + ({n_1} - 1)S^{'}(u_2|{G_2}).\]
\end{lemma}

\begin{definition}
Let $G$ be a bicyclic graph $G$ with $V(G) = \{u, v,v_1, \dots, v_s\} \cup V(C_p) \cup V(C_q)$, for which
  $v$ is a vertex of degree $s+1$
  such that $vv_1, vv_2, \ldots, vv_s$ are pendent edges incident with $v$,
   and $u$ is the neighbor of $v$ distinct from $v_i$ that is on the cycle $C_q$.
   The other cycle $C_p$ only has one common vertex $w$ with
   $C_q$.
   We form a graph ${G^{'}} = \sigma (G,v)$ by deleting the edges $vv_1,vv_2,\ldots,vv_s$
   and adding new edges $uv_1,uv_2,\ldots,uv_s$.
   We say that ${G^{'}}$ is a $\sigma$-transform of the graph $G$,
   see Fig. 3.
\end{definition}

\begin{minipage}[t]{0.5\linewidth}
  \centering
    \includegraphics[scale=0.75]{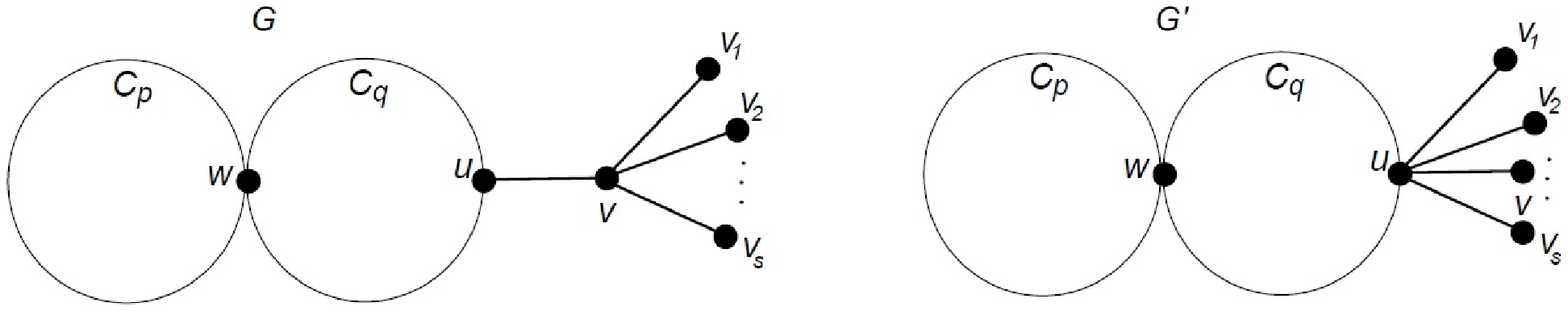}
  \vspace{-2em} \centerline{ ~~~~~~~~~~~~~~~~~~~~~~~~~~~~~~~~~~~~~~~~~~~~~~~~~~~~~~~~~Fig. 3. $\sigma$-transform of a vertex $v.$}
\end{minipage}%
\vspace{1em}

\begin{lemma}
Consider the graphs defined in Definition 2.1 and
 let ${G^{'}} = \sigma (G,v)$ be a $\sigma$-transform of the bicyclic graph $G$, for convenience, let $r(w,u)=t$.
  Then ${D_R}(G) > {D_R}({G^{'}})$.
\end{lemma}

\noindent \textbf{Proof}.
 Let $H$ be a subgraph $G[\{u, v, v_1, \dots, v_s\} \cup V(C_q)]$ and $T = G[\{u,v, v_1, \dots,v_s\}]$
 ($H'$ be a subgraph $G'[\{u, v, v_1, \dots, v_s\} \cup V(C_q)]$ and $T' = G'[\{u,v, v_1, \dots,v_s\}]$, respectively).
  Note that $H$ and $C_p$ ($H'$ and $C_p$, respectively) share a common vertex $w$.
 By Lemmas 2.2 and 2.3, we get
$${D_R} ({G} )={D_R} ({C_p} )+{D_R} ({H}
)+2({s+1+q})r({w|C_p })+{2p}r({w|H
})+({q+s})S^{'}({w|C_p})+({p-1})S^{'}({w|H}),$$
$${D_R} ({H} )={D_R} ({C_q} )+{D_R} ({T} )+2({s+1})r({u|C_q
})+{2q}r({u|T })+({s+1})S^{'}({u|C_q})+({q-1})S^{'}({u|T}),$$
$${D_R} ({G}^{'} )={D_R} ({C_p} )+{D_R} ({H}^{'}
)+2({s+1+q})r({w|C_p })+{2p}r({w|H^{'}
})+({q+s})S^{'}({w|C_p})+({p-1})S^{'}({w|H^{'}}),$$
$$ {D_R}
({H}^{'} )={D_R} ({C_q} )+{D_R} (T')+2({s+1})r({u|C_q })+{2q}r({u|T'
})+({s+1})S^{'}({u|C_q})+({q-1})S^{'}({u|T'}),$$
\begin{equation*}
\begin{split}
{D_R} ({G} )-{D_R} ({G}^{'} )&={2q}r({u|T
})+({q-1})S^{'}({u|T})+{2p}r({w|H })+({p-1})S^{'}({w|H})\\
&\quad-{2q}r({u|T'
})-({q-1})S^{'}({u|T'})-{2p}r({w|H^{'}})-({p-1})S^{'}({w|H^{'}})\\
&=2q(2s+1)+(q-1)(s+1+2s)+2p\Big[\frac{{q^2-1}}{6}+t+1+(t+2)s\Big]\\
&\quad+(p-1)\Big[\frac{q^2-1}{3}+(t+1)(s+1)+(t+2)s\Big]-2q(s+1)-(q-1)(s+1)\\
&\quad-2p\Big[\frac{q^2-1}{6}+(t+1)(s+1)\Big]-(p-1)\Big[\frac{q^2-1}{3}+(t+1)(s+1)\Big]\\
&=2ps+(p-1)(t+2)s+2s(2q-1)>0.
\end{split}
\end{equation*}

\begin{lemma} Let $G_0$ be a bicyclic graph with $V(G) = \{v_1,\dots,v_s\} \cup V(C_p) \cup V(C_q)$, for which
$u$ is a vertex of degree $s$ in the cycle $C_q$ of the
bicyclic graph $G_0$, and $uv_1,uv_2,...,uv_s$ are pendent edges
incident with $u$. The other cycle $C_p$ only has one common
vertex $w$ with $C_q$. Let graph $G_1$ delete the edges
$uv_1,uv_2,...,uv_s$, and add new edges $wv_1,wv_2,...,wv_s$. For
convenience, let $r(w,u)=t$. Then $D_R(G_0)> D_R(G_1)$.
\end{lemma}

\noindent  \textbf{Proof}.
 Let $H_0$ be a subgraph $G[\{ v_1, \dots, v_s\} \cup V(C_q)]$ and $S = G[\{u, v_1, \dots,v_s\}]$ ($H_1$ be a subgraph $G_1[\{  v_1, \dots, v_s\} \cup V(C_q)]$ and $S_1 = G'[\{w,v, v_1, \dots,v_s\}]$, respectively).
By Lemmas 2.2 and 2.3, we get
$${D_R} ({G_0} )={D_R} ({C_p} )+{D_R} ({H_0}
)+2({s+q})r({w|C_p})+{2p}r({w|H_0
})+({q+s-1})S^{'}({w|C_p})+({p-1})S^{'}({w|H_0}),$$
$${D_R} ({H_0} )={D_R} ({C_q} )+{D_R} ({S} )+{2s}r({u|C_q
})+{2q}r({u|S })+{s}S^{'}({u|C_q})+({q-1})S^{'}({u|S}),$$
$${D_R} ({G_1} )={D_R} ({C_p} )+{D_R} ({H_1} )+2({s+q})r({w|C_p
})+{2p}r({w|H_1})+({q+s-1})S^{'}({w|C_p})+({p-1})S^{'}({w|H_1}),$$
$${D_R} ({H_1} )={D_R} ({C_q} )+{D_R} ({S_1} )+2{s}r({w|C_q
})+{2q}r({w|S_1})+{s}S^{'}({w|C_q})+({q-1})S^{'}({w|S_1}),$$
\begin{equation*}
\begin{split}
{D_R} ({G_0} )-{D_R} ({G_1} )&={2p}r({w|H_0
})+({p-1})S^{'}({w|H_0})-{2p}r({w|H_1 })-({p-1})S^{'}({w|H_1})\\
&=2p\Big[\frac{q^2-1}{6}+(t+1)s\Big]+(p-1)\Big[\frac{q^2-1}{3}+(t+1)s\Big]\\
&\quad-2p\Big[\frac{q^2-1}{6}+s\Big]-(p-1)\Big[\frac{q^2-1}{3}+s\Big]\\
&=(3p-1)ts>0.
\end{split}
\end{equation*}

\begin{definition}
 Let $G$ be a bicyclic graph $G$ with $V(G) = \{u, v,v_1, \dots, v_s\} \cup V(C_p) \cup V(C_q)$, for which
  $v$ is a vertex of degree $s+1$
  such that $vv_1, vv_2, \ldots, vv_s$ are pendent edges incident with $v$,
   and $u$ is the neighbor of $v$ distinct from $v_i$ that is on the cycle $C_q$.
   The other cycle $C_p$ only has one common vertex $w$ with
   $C_q$.
   We form a graph ${G^{''}} = \pi (G,v)$ by deleting the edges $vv_1, vv_2,\ldots, vv_s$
   and connecting $v_i$ and all the isolated vertices into a path $vv_1 \dots v_s$.
   We say that ${G^{''}}$ is a $\pi$-transform of the graph $G$, see Fig. 4.
\end{definition}

\begin{minipage}[t]{0.5\linewidth}
  \centering
    \includegraphics[scale=0.75]{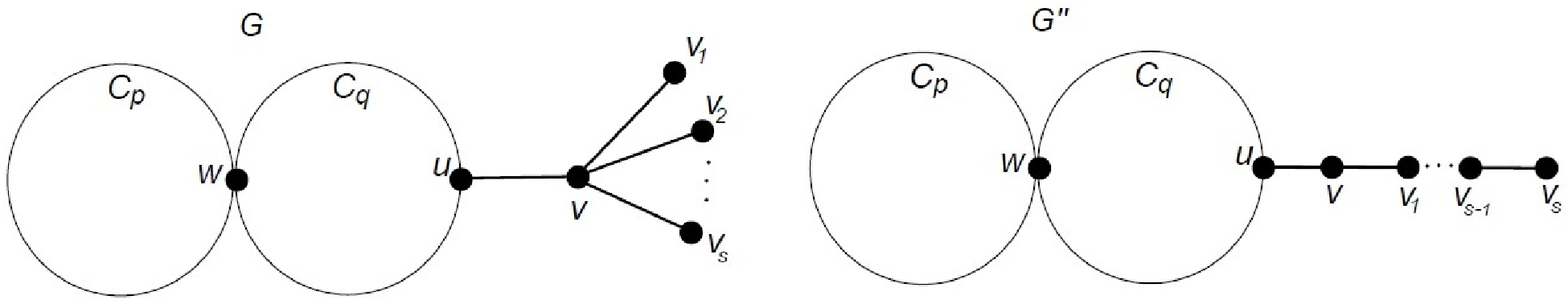}
  \vspace{-2em} \centerline{ ~~~~~~~~~~~~~~~~~~~~~~~~~~~~~~~~~~~~~~~~~~~~~~~~~~~~~~~~~Fig. 4. $\pi$-transform of a vertex $v.$}
\end{minipage}%
\vspace{1em}

\begin{lemma}
Consider the graphs defined in Definition 2.1 and
let ${G^{''}} = \pi (G,v)$ be a $\pi$-transform of the bicyclic graph G, for convenience, let $r(w,u)=t$.
  Then ${D_R}(G) < {D_R}({G^{''}})$.
\end{lemma}

\noindent \textbf{Proof}.
 Let $H$ be a subgraph $G[\{u, v, v_1, \dots, v_s\} \cup V(C_q)]$ and $T = G[\{u,v, v_1, \dots,v_s\}]$ ($H''$ be a subgraph $G''[\{u, v, v_1, \dots, v_s\} \cup V(C_q)]$ and $P = G'[\{u,v, v_1, \dots,v_s\}]$, respectively).
By Lemmas 2.2 and 2.3, we get
$${D_R} ({G} )={D_R} ({C_p} )+{D_R} ({H}
)+2({s+1+q})r({w|C_p })+{2p}r({w|H
})+({q+s})S^{'}({w|C_p})+({p-1})S^{'}({w|H}),$$
$${D_R} ({H} )={D_R} ({C_q} )+{D_R} ({T} )+2({s+1})r({u|C_q
})+{2q}r({u|T })+({s+1})S^{'}({u|C_q})+({q-1})S^{'}({u|T}),$$
$${D_R} ({G}^{''} )={D_R} ({C_p} )+{D_R} ({H}^{''}
)+2({s+1+q})r({w|C_p })+{2p}r({w|H^{''}
})+({q+s})S^{'}({w|C_p})+({p-1})S^{'}({w|H^{''}}),$$
$$ {D_R}
({H}^{''} )={D_R} ({C_q} )+{D_R} (P)+2({s+1})r({u|C_q
})+{2q}r({u|P})+({s+1})S^{'}({u|C_q})+({q-1})S^{'}({u|P}),$$
\begin{equation*}
\begin{split}
{D_R} ({G} )-{D_R} ({G}^{''} )&={D_R} ({T} )+{2q}r({u|T
})+({q-1})S^{'}({u|T})+{2p}r({w|H })+({p-1})S^{'}({w|H})\\
&\quad-{D_R} ({P} )-{2q}r({u|P
})-({q-1})S^{'}({u|P})-{2p}r({w|H^{''}})-({p-1})S^{'}({w|H^{''}})\\
&=(3s^2-s)+2q(2s+1)+(q-1)(s+1+2s)+2p\Big[\frac{{q^2-1}}{6}+t+1+(t+2)s\Big]\\
&\quad+(p-1)\Big[\frac{q^2-1}{3}+(t+1)(s+1)+(t+2)s\Big]-(\frac{2}{3}s^3+s^2+\frac{1}{3})\\
&\quad-2q\frac{(1+s+1)(s+1)}{2}-(q-1)\Big[(s+1)s+s+1\Big]\\
&\quad-2p\Big[\frac{q^2-1}{6}+t(s+1)+\frac{(1+s+1)(s+1)}{2}\Big]\\
&\quad-(p-1)\Big[\frac{q^2-1}{3}+(t+1+t+s)s+t+s+1\Big]\\
&=(s-s^2)(2q+2p-2)-\frac{2}{3}s(s-1)(s-2)<0.
\end{split}
\end{equation*}

\begin{lemma}
Let $G_0$ be a bicyclic graph with the vertex set $V(C_p) \cup
V(C_q) \cup V(P_{s+1})$, in which $V(C_p) \cap V(P_{s+1}) = \{v\}$
and $V(C_q) \cap V(P_{s+1}) = \{w\}$. For $wa \in E(P_{s+1})$ and
$u \in V(C_q)$, let $G_1 = (G - \{aw\}) \cup \{ua\}$. For
convenience, let $r(u,w)=t$. Then $D_R(G_0)> D_R(G_1). $
\end{lemma}

\noindent  \textbf{Proof}.  By Lemmas 2.2 and 2.3, we get
$${D_R} ({G_0}
)={D_R} ({H_0} )+{D_R} ({H} )+{2(p+s-1)}r({a|H_0 })+2(q+1)r({a|H
})+({p+s-2})S^{'}({a|H_0})+{q}S^{'}({a|H}),$$
$${D_R} ({G}_1 )={D_R} ({H_1} )+{D_R} ({H} )+{2(p+s-1)}r({w|H_1
})+2(q+1)r({w|H })+({p+s-2})S^{'}({w|H_1})+{q}S^{'}({w|H}),$$
\begin{equation*}
\begin{split}
{D_R} ({G_0} )-{D_R} ({G}_1 )&=2(p+s-1)r({a|H_0
})+({p+s-2})S^{'}({a|H_0})\\
&\quad-2(p+s-1)r({w|H_1})-({p+s-2})S^{'}({w|H_1})\\
&=2(p+s-1)(1+\frac{q^2-1}{6}+q)+(p+s-2)(\frac{q^2-1}{3}+3+2q)\\
&\quad-2(p+s-1)(\frac{q^2-1}{6}+t+1)-(p+s-2)(\frac{q^2-1}{3}+2t+1)\\
&=2(p+s-1)(q-t)+(p+s-2)(2q-2t+2)>0.
\end{split}
\end{equation*}

\section{Main results}

In this section, we will characterize $n$-vertex bicyclic graphs
 with exactly two cycles having minimum and maximum degree resistance
distances.

\begin{theorem}
 Let ${G} \in \mathscr {B}_n^{p,q} $and  $G\neq S_n^{p,q}.$ Then
 $D_R(G)>D_R(S_n^{p,q}).$
\end{theorem}

\noindent  \textbf{Proof}. Suppose that a bicyclic graph $G_0$ has
minimal degree resistance distance among graphs in $\mathscr
{B}_n^{p,q}$ . For $G_0$, we prove the following claims.

\noindent \textbf{Claim 1}. In Fig. 1, $T_{v_i},T_{u_j}$ and
$T_{w_k}$ are all stars with their centers at $v_i,u_j$ and $w_k$
for each $i,j$ and $k$.

Without loss of generality, suppose that tree $T_{v_i}$ is not a
star. Let $G_1$ be constructed from $G_0$ by deleting all the
edges of $T_{v_i}$ and connecting all the isolated vertices to
$v_i$; that is, $T_{v_i}$ is a star in $G_1$ with its center at
$v_i$ and denote it by $S_{v_i}$. By Lemma 2.4,
$D_R(G_0)>D_R(G_1)$, which contradicts the choice of $G_0$. Hence
Claim 1 holds.
 \hfill$\blacksquare$

\noindent \textbf{Claim 2}. The length of $P$ is 0.

Suppose to the contrary that the length of $P$ is $k(k\geq1)$.
Assume that $v_1=w_0,u_1=w_k$. Let $e=w_iw_{i+1}$ be an edge of
$P$. Let $G_2$ be the graph obtained from $G_0$ by first
contracting $e$ and then attaching a pendent edge $w_ia$ to $w_i$.
Assume that $G_{01}$ and $G_{02}$ are two components of $G_0-e$
and $G_{21}$ and $G_{22}$ are copies of $G_{01}$ and $G_{02}$ in
$G_2$, respectively. See Fig. 5.

\begin{figure}[ht]
\centering
  \includegraphics[width=\textwidth]{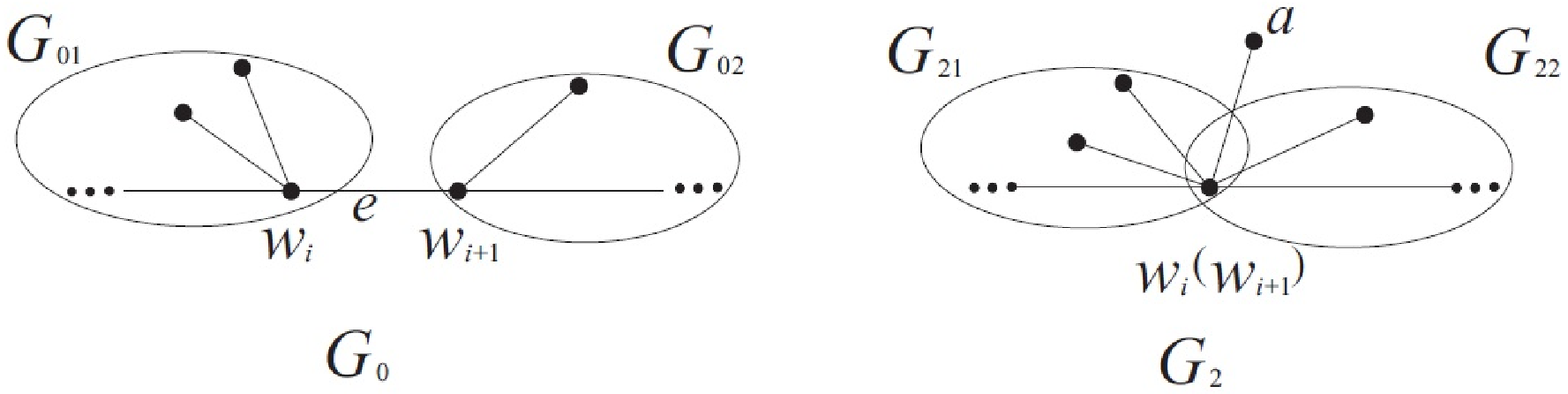}
 \vspace{-2em} \centerline{ ~~~Fig. 5. Graphs $G_{0}$ and $G_{2}$.}
\end{figure}

\vspace{2em}

In the following, we prove $D_R(G_2)<D_R(G_0)$.

\noindent \textbf{Proof}. By Lemma 2.3, we get
\begin{equation*}
\begin{split}
{D_R} ({G}_0 )&={D_R} ({G}_{01} )+{D_R}
(G_{02}+w_iw_{i+1})+2({m_{02}+1})r({w_i|G_{01}})+{2m_{01}}r({w_i|G_{02}+w_iw_{i+1}
})\\
&\quad+n_{02}S^{'}({w_i|G_{01}})+({n_{01}-1})S^{'}({w_i|G_{02}+w_iw_{i+1}}),
\end{split}
\end{equation*}
\begin{equation*}
\begin{split}
{D_R} ({G}_2 )&={D_R} ({G}_{21} )+{D_R}
(G_{22}+w_ia)+2({m_{22}+1})r({w_i|G_{21}})+{2m_{21}}r({w_i|G_{22}+w_ia
})\\
&\quad+n_{22}S^{'}({w_i|G_{21}})+({n_{21}-1})S^{'}({w_i|G_{22}+w_ia}),
\end{split}
\end{equation*}
\begin{equation*}
\begin{split}
{D_R} ({G}_0 )-{D_R} ({G}_2 )&={2m_{01}}r({w_i|G_{02}+w_iw_{i+1}
})+({n_{01}-1})S^{'}({w_i|G_{02}+w_iw_{i+1}})\\
&\quad-{2m_{21}}r({w_i|G_{22}+w_ia
})-({n_{21}-1})S^{'}({w_i|G_{22}+w_ia})\\
&={2m_{01}}\Big[r({w_{i+1}|G_{02})+n+1}\Big]+({n_{01}-1})\Big[S^{'}({w_{i+1}|G_{02})+n+n+1}\Big]\\
&\quad-{2m_{21}}\Big[r({w_{i+1}|G_{22})+1
}\Big]-({n_{21}-1})\Big[S^{'}({w_{i+1}|G_{22}})+n+1\Big]\\
&=(2m+n-1)n>0.
\end{split}
\end{equation*}

We obtain $D_R(G_2)<D_R(G_0)$. This contradicts the hypothesis.
Hence Claim 2 holds.
 \hfill$\blacksquare$

\noindent \textbf{Claim 3}. In Fig. 1, if $p+q\leq n$, then only
$T_{v_1}(T_{v_1}=T_{u_1})$ is nontrivial.

Without loss of generality, suppose to the contrary that tree
$T_{u_i}(i\neq1)$ is nontrivial. By Lemma 2.5,
$D_R(G_0)>D_R(G_1)$, which contradicts the choice of $G_0$. Hence
Claim 3 holds.

Claims 1-3 yield Theorem 3.1.
 \hfill$\blacksquare$

 \begin{figure}[ht]
\centering
  \includegraphics[width=\textwidth]{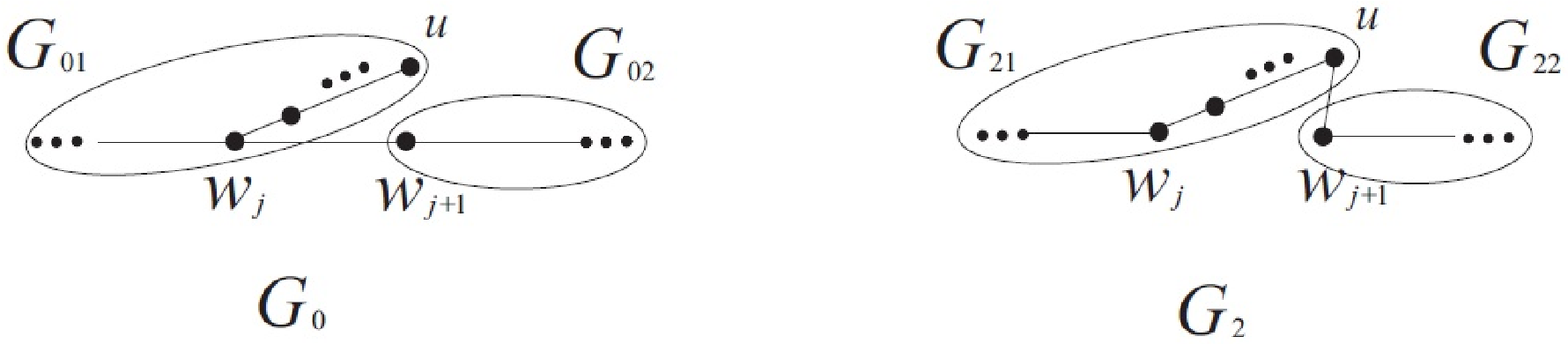}
 \vspace{-2em} \centerline{ ~~~Fig. 6. Graphs $G_{0}$ and $G_{2}$.}
\end{figure}

\vspace{2em}

\begin{theorem}.
 Let ${G} \in \mathscr {B}_n^{p,q} $and  $G\neq P_n^{p,q}.$ Then
 $D_R(G)<D_R(P_n^{p,q}).$
\end{theorem}

\noindent \textbf{Proof}. Suppose that a bicyclic graph $G_0$ has
maximal degree resistance distance among graphs in $\mathscr
{B}_n^{p,q}$. For $G_0$, we prove the following Claims.

\noindent \textbf{Claim 1}. In Fig. 1, $T_{v_i},T_{u_j}$ and
$T_{w_k}$ are all paths with their end vertices $v_i,u_j$ and
$w_k$ for each $i,j$ and $k$.

Without loss of generality, suppose that tree $T_{v_i}$ is not a
path. Let $G_1$ be the graph constructed from $G_0$ by deleting
all the edges of $T_{v_i}$ and connecting $v_i$ and all the
isolated vertices into a path; that is, $T_{v_i}$ is a path with
end vertex $v_i$ in $G_1$ and denote it by $P_{v_i}$. By Lemma
2.6, $D_R(G_1)>D_R(G_0)$, which contradicts the choice of $G_0$.
Hence Claim 1 holds.

\noindent \textbf{Claim 2}. Assume that $T_{w_0}=T_{v_1}$ and
$T_{w_m}=T_{u_1}$, then $T_{w_i}$ is trivial ($0\leq i \leq m$).

  If not, without loss of generality, suppose that there is
  nontrivial $T_{w_j}$. By Claim 1, we know that $T_{w_j}$ is a path
  with $w_j$ as its end vertex and assume that $u$ is the other
  end vertex. Let $G_2=G_0-w_j{w_{j+1}}+u{w_{j+1}}$(if
  $j=m$, $G_2=G_0-w_{j-1}w_j+uw_{j-1}$). Assume that $G_{01}$ and $G_{02}$
  are two components of $G_0-w_j{w_{j+1}}$ and $G_{21}$ and $G_{22}$ are
  two components of $G_2-uw_{j+1}$. See Fig. 4.

In the following, we prove $D_R(G_2)>D_R(G_0)$.

\noindent \textbf{Proof}. By Lemma 2.3, we get
\begin{equation*}
\begin{split}
{D_R} ({G}_0 )&={D_R} (G_{01}+w_jw_{j+1})+{D_R} ({G}_{02}
)+{2m_{02}}r({w_{j+1}|G_{01}+w_jw_{j+1}
})+2({m_{01}+1})r({w_{j+1}|G_{02}})\\
&\quad+({n_{02}-1})S^{'}({w_{j+1}|G_{01}+w_jw_{j+1}})+n_{01}S^{'}({w_{j+1}|G_{02}}),
\end{split}
\end{equation*}
\begin{equation*}
\begin{split}
{D_R} ({G}_2 )&={D_R} (G_{21}+uw_{j+1})+{D_R} ({G}_{22}
)+{2m_{22}}r({w_{j+1}|G_{21}+uw_{j+1}
})+2({m_{21}+1})r({w_{j+1}|G_{22}})\\
&\quad+({n_{22}-1})S^{'}({w_{j+1}|G_{21}+uw_{j+1}})+n_{21}S^{'}({w_{j+1}|G_{22}}),
\end{split}
\end{equation*}
\begin{equation*}
\begin{split}
{D_R} (G_{01}+w_jw_{j+1})&={D_R} ({G}_{01} )+{D_R} (w_jw_{j+1})
+2r({w_{j}|G_{01}})+2m_{01}r(w_j|w_jw_{j+1})\\
&\quad+S^{'}(w_{j}|G_{01})+(n_{01}-1)S^{'}(w_{j}|w_jw_{j+1}),
\end{split}
\end{equation*}
\begin{equation*}
\begin{split}
{D_R} (G_{21}+uw_{j+1})&={D_R} ({G}_{21} )+{D_R} (uw_{j+1})
+2r({u}|G_{21})+2m_{21}r(u|uw_{j+1})\\
&\quad+S^{'}(u|G_{21})+(n_{21}-1)S^{'}(u|uw_{j+1}),
\end{split}
\end{equation*}
\begin{equation*}
\begin{split}
{D_R} ({G}_0 )-{D_R} ({G}_2 )&=
2r(w_j|{G_{01}})+S^{'}(w_j|{G_{01}})+{2m_{02}}r({w_{j+1}|G_{01}+w_jw_{j+1}
})\\
&\quad+({n_{02}-1})S^{'}({w_{j+1}|G_{01}+w_jw_{j+1}})-2r(u|{G_{21}})-S^{'}(u|{G_{21}})\\
&\quad-{2m_{22}}r({w_{j+1}|G_{21}+uw_{j+1}
})-({n_{22}-1})S^{'}({w_{j+1}|G_{21}+uw_{j+1}})\\
&<{2m_{02}}\Big[r({w_{j+1}|G_{01})+1}\Big]+({n_{02}-1})\Big[S^{'}({w_{j+1}|G_{01})+3}\Big]\\
&\quad-{2m_{22}}\Big[r({w_{j+1}|G_{21})+1
}\Big]-({n_{22}-1})\Big[S^{'}({w_{j+1}|G_{21}})+2\Big]\\
&=(n-1)\Big[S^{'}({w_{j+1}|G_{01})})-S^{'}({w_{j+1}|G_{21}})+1\Big]<0.
\end{split}
\end{equation*}

We obtain $D_R(G_2)>D_R(G_0)$. This contradicts the hypothesis.
Hence Claim 2 holds.

\noindent \textbf{Claim 3}. In Fig. 1, if $p+q\leq n$, then
$T_{v_i}$ and $T_{u_j}$ are trivial for each $i$ and $j$.

Without loss of generality, suppose to the contrary that
$T_{v_i}(i\neq1)$ is nontrivial. By Lemma 2.7,
$D_R(G_0)<D_R(G_1)$, which contradicts the choice of $G_0$. Hence
Claim 3 holds.

Claims 1-3 yield Theorem 3.2.
 \hfill$\blacksquare$

We now compute $D_R(S_n^{p,q})$ and $D_R(P_n^{p,q})$.

For $S_n^{p,q}$, see Fig. 2(a).
$${D_R} ({S_n^{p,q}} )={D_R} ({C_q} )+{D_R} ({H}
)+2{(p+s)}r({v_1|C_q })+2qr({v_1|H
})+({p+s-1})S^{'}({v_1|C_q})+({q-1})S^{'}({v_1|H}),$$
$${D_R} ({H} )={D_R} ({C_p} )+{D_R} ({S} )+2sr({v_1|C_p
})+{2p}r({v_1|S })+sS^{'}({v_1|C_p})+({p-1})S^{'}({v_1|S}),$$
$$r(v_1|H)=r(v_1|C_p)+r(v_1|S),$$
$$S^{'}({v_1|H})=S^{'}({v_1|C_p})+S^{'}({v_1|S}),$$
\begin{equation*}
\begin{split}
{D_R} ({S_n^{p,q}} )&=\frac{q^3-q}{3}
+\frac{p^3-p}{3}+3(s+1)^2-7(s+1)+4+2s\frac{p^2-1}{6}+2ps\\
&\quad+s\frac{p^2-1}{3}+(p-1)s+2(p+s)\frac{q^2-1}{6}+2q\frac{p^2-1}{6}+2qs\\
&\quad+(p+s-1)\frac{q^2-1}{3}+(q-1)\frac{p^2-1}{3}+(q-1)s\\
&=\frac{1}{3}(q^3+p^3+2qp^2+2q^2p+2s{p}^2+2sq^2+9s{p}+9sq-q^2-p^2-3q-3p+9s^2-13s+2).\\
\end{split}
\end{equation*}
And $s=n+1-p-q$, so we get
\begin{equation}
\begin{split}
{D_R} ({S_n^{p,q}}
)=\frac{1}{3}[-p^3-q^3+(2n+1)(p^2+q^2)+(1-9n)(p+q)+9n^2+5n-2].
\end{split}
\end{equation}

For $P_n^{p,q}$, see Fig. 2(b).
\begin{equation*}
\begin{split}
{D_R} (P_n^{p,q} )&={D_R} ({C_p} )+{D_R} ({H} )+{2(q+m)}r({v_1|C_p
})+2pr({v_1|H
})+({q+m-1})S^{'}({v_1|C_p})+({p-1})S^{'}({v_1|H}),\\
{D_R} ({H} )&={D_R} ({C_q} )+{D_R} ({P_{m+1}} )+2qr({u_1|P_{m+1}
})+{2m}r({u_1|C_q })+(q-1)S^{'}({u_1|P_{m+1}})+sS^{'}({u_1|C_q}),\\
\end{split}
\end{equation*}
$$r(v_1|H)=r(v_1|C_q)+r(v_1|P_{m+1}),$$
$$S^{'}({v_1|H})=S^{'}({v_1|C_q})+S^{'}({v_1|P_{m+1}}),$$
\begin{equation*}
\begin{split}
{D_R} ({P_n^{p,q}} )&=\frac{q^3-q}{3}
+\frac{p^3-p}{3}+\frac{2}{3}m^3+m^2+\frac{1}{3}m+2q\frac{(1+m)m}{2}+2m\frac{q^2-1}{6}+(q-1)m^2+\frac{q^2-1}{3}m\\
&\quad+2(q+m)\frac{p^2-1}{6}+2p[\frac{(1+m)m}{2}+m(q-1)+\frac{q^2-1}{6}]\\
&\quad+(q+m-1)\frac{p^2-1}{3}+(p-1)[(1+m-1)(m-1)+3m+2m(q-1)+\frac{q^2-1}{3}]\\
&=\frac{1}{3}\Big[p^3+q^3+(2q+2m-1)p^2+(2p+2m-1)q^2+(6m^2-3m-3)(p+q)\\
&\quad+12mpq+2m^3-3m^2-3m+2\Big].\\
\end{split}
\end{equation*}
And $m=n+1-p-q$, so we get
\begin{equation}
\begin{split}
{D_R} ({P_n^{p,q}}
)&=\frac{1}{3}\Big[3p^3+3q^3+(4n+5)(p^2+q^2)+(3n+3)(p+q)+2n^3+3n^2-3n-2\Big].
\end{split}
\end{equation}
 \hfill$\blacksquare$

\section{Bicyclic graphs with extremal degree resistance distance}

  By Theorems 3.1 and 3.2, $n$-vertex bicyclic graphs in $\mathscr {B}_n^{p,q}$
with minimal and maximal degree resistance distance must belong to
the classes of $S_n^{p,q}$ and $P_n^{p,q}$, respectively. In what
follows, we will determine those which has the extremal degree
resistance distance among graphs in $\mathscr {B}_n^{p,q}$.

\begin{theorem}
 Among all $n$-vertex bicyclic graphs, the graph $S_n^{3,3}$ has the minimum degree resistance distance.
 \begin{center}
 $D_R(S_n^{3,3})=3n^2-\frac{13}{3}n-\frac{32}{3},~~~~~~~n\geq5$.
\end{center}
\end{theorem}

\noindent \textbf{Proof}. Let $u_1,u_2, w$ be three successive
vertices lying on the cycle $C_p$ of the bicyclic graph $G_1$. The
other cycle $C_q$ only has one common vertex $w$ with $C_p$. And
$wv_1, wv_2,..., wv_s$ are pendent edges incident with $w$. Let
the graph $G_2$ is obtained by deleting the edges $u_1u_2$ and
adding the edge $wu_2$. Then $D_R(G_2)<D_R(G_1)$.
$${D_R} ({G_1} )={D_R} ({C_p} )+{D_R}
({H_1} )+{2p}r({w|H_1 })+2({s+q})r({w|C_p
})+({p-1})S^{'}({w|H_1})+({q+s-1})S^{'}({w|C_p}),$$
$${D_R} ({G_2}
)={D_R} ({C_p^{'}} )+{D_R} ({H_1} )+{2p}r({w|H_1
})+2({s+q})r({w|C_p^{'}
})+({p-1})S^{'}({w|H_1})+({q+s-1})S^{'}({w|C_p^{'}}),$$
\begin{equation*}
\begin{split}
{D_R} ({G_2} )-{D_R} ({G_1} )&={D_R} ({C_p^{'}}
)+2({s+q})r({w|C_p^{'} })+({q+s-1})S^{'}({w|C_p^{'}})\\
&\quad-{D_R} ({C_p} )-2({s+q})r({w|C_p
})-({q+s-1})S^{'}({w|C_p})\\
&={D_R}({C_{p-1}})+{D_R}({u_1w})+2(p-1)r(w|u_1w)+2r(w|C_{p-1})\\
&\quad+(p-1)S^{'}({w|u_1w})+S^{'}({w|C_{p-1}})-{D_R} ({C_p}
)+2(s+q)r(w|C_p^{'})\\
&\quad-2(s+q)r(w|C_{p})+(s+q-1)S^{'}({w|C_p^{'}})-(s+q-1)S^{'}({w|C_p})\\
&=\frac{(p-1)^3-(p-1)}{3}+2+2p-2+2\frac{(p-1)^2-1}{6}+p-1\\
&\quad+\frac{(p-1)^2-1}{3}-\frac{p^3-p}{3}+2(s+q)\Big[\frac{(p-1)^2-1}{6}+1-\frac{p^2-1}{6}\Big]\\
&\quad+(s+q-1)\Big[\frac{(p-1)^2-1}{3}+1-\frac{p^2-1}{3}\Big]\\
&=-\frac{1}{3}p^2-\frac{8}{3}p-1-\frac{4}{3}+\frac{2}{3}p+\frac{s}{3}(11-4p)+\frac{q}{3}(11-4p)\\
&=-\frac{1}{3}p^2-2p-\frac{7}{3}+\frac{s}{3}(11-4p)+\frac{q}{3}(11-4p)<0.
\end{split}
\end{equation*}
Combining the above discussions, Eq. (1) and $p=q=3$, we can get\\
$$D_R(S_n^{3,3})=3n^2-\frac{13}{3}n-\frac{32}{3}, ~~~~~~~~n\geq5.$$
 \hfill$\blacksquare$

\begin{theorem}
 Among all $n$-vertex bicyclic graphs, the graph $P_n^{3,3}$ has the maximal degree resistance distance.
 $$D_R(P_n^{3,3})=\frac{2}{3}n^3+n^2-19n+\frac{88}{3},~~~~~~~��n\geq5��.$$
\end{theorem}

\noindent  \textbf{Proof}. Let $u_1, w, u_2$ be three successive
vertices lying on the cycle $C_p$ of the bicyclic graph $G_3$. The
cycles $C_p$ and $C_q$ are linked with two end vertexes $v$ and
$w$ of $P_{s+1}$. Let the graph $G_4$ is obtained by deleting the
edge $wu_2$ and adding the edge $u_1u_2$. Then
$D_R(G_4)>D_R(G_3)$.
$${D_R} ({G}_3)={D_R} ({C_p} )+{D_R} ({H}
)+{2(q+s)}r({w|C_p })+2pr({w|H
})+({q+s-1})S^{'}({w|C_p})+({p-1})S^{'}({w|H}),$$
$${D_R}
({G}_4)={D_R} ({C_p^{'}} )+{D_R} ({H} )+{2(q+s)}r({w|C_p^{'}
})+2pr({w|H })+({q+s-1})S^{'}({w|C_p^{'}})+({p-1})S^{'}({w|H}),$$
\begin{equation*}
\begin{split}
{D_R} ({C_p^{'}} )&={D_R} ({C_{p-1}} )+{D_R} ({wu_1}
)+{2}r({u_1|C_p^{'} })+2(p-1)r({u_1|wu_1
})+S^{'}({u_1|C_p^{'}})+({p-2})S^{'}({u_1|wu_1})\\
&=\frac{(p-1)^3-(p-1)}{3}+2+2\frac{(p-1)^2-1}{6}+2(p-1)+\frac{(p-1)^2-1}{3}+(p-2)\\
&=\frac{p^3}{3}+\frac{p^2}{3}+\frac{8p}{3}+\frac{7}{3},
\end{split}
\end{equation*}
$$r({w|C_p^{'} })=1+\frac{(p-1)^2-1}{3}+(p-1)
=\frac{p^2}{6}+\frac{2p}{3},$$
$$S^{'}({w|C_p^{'}})=3+
\frac{(p-1)^2-1}{3}+2(p-1)=\frac{p^2-2p}{3}+2p+1,$$
\begin{equation*}
\begin{split}
{D_R} ({G_4} )-{D_R} ({G}_3)&={D_R} ({C_p^{'}}
)+2(q+s)r({w|C_p^{'} })+(q+s-1)S^{'}(w|C_p^{'})\\
&\quad-{D_R} ({C_p})-2(q+s)r({w|C_p})-(q+s-1)S^{'}(w|C_p)\\
&=\frac{p^3}{3}+\frac{p^2}{3}+\frac{8p}{3}-\frac{7}{3}+2(q+s)(\frac{p^2}{6}+\frac{2p}{3})\\
&\quad+(q+s-1)(\frac{p^2-2p}{3}+2p+1)-\frac{p^3-p}{3}\\
&\quad-2(q+s)\frac{p^2-1}{6}-(q+s-1)\frac{p^2-1}{3}\\
&=\frac{p^3}{3}+3p+2(q+s)(\frac{2p}{3}+\frac{1}{6})+(q+s-1)(\frac{4p}{3}+\frac{4}{3})-\frac{7}{3}>0.
\end{split}
\end{equation*}

Combining above discussion, and by Eq. (2) and $p=q=3$. So we can
get
$$D_R(P_n^{3,3})=\frac{2}{3}n^3+n^2-19n+\frac{88}{3},~~~~~~~��n\geq5��.$$
 \hfill$\blacksquare$

\section{Concluding remarks}

 In this paper, we completely characterized the
bicyclic graphs with exactly two cycles having extremal degree
resistance distances. Although our research is restricted to a
small family of graphs, it motivates one to further extend this to
the larger classes of graphs. For example, can one determine the
maximum or minimum degree resistance distance for tricyclic graphs
or any general graphs? We plan to investigate this question for
larger classes of graphs for future research in this area.

\section*{Acknowledgments}

The work of J.B. Liu was partly supported by the Natural Science
Foundation of Anhui province of China under grant No. KJ2013B105,
the Natural Science Foundation for the Higher Education
Institutions of Anhui province of China under grant No.
KJ2015A331, and the National Science Foundation of China under
grant No's. 11471016 and 11401004. The work of X. F. Pan was
partly supported by the National Science Foundation of China under
grant No's. 10901001, 11171097, and 11371028. The work of S. Hayat
is supported by CAS-TWAS president's fellowship at USTC, China.

\end{document}